\documentclass[12pt]{amsart}

\usepackage{epsfig}
\usepackage{graphics}
\usepackage{latexsym}
\usepackage{psfrag}
\usepackage{amssymb}
\usepackage{amsmath}
\usepackage{amsthm} 

\renewcommand{\setminus}{{\smallsetminus}}

\newcommand{\BB}{{\mathbb{B}}}

\newcommand{\pt}{{\rm{pt}}} 
\newcommand{\bdy}{{\partial}} 

\theoremstyle{plain}
\newtheorem{theorem}{Theorem}[section]

\theoremstyle{definition}
\newtheorem*{define}{Definition}
\newtheorem*{claim}{Claim}

\newsavebox{\savepar}

\begin{document}

\title{Thin Position for Tangles} 
\date{December 1, 2001}
\address{Mathematics Department, University of Illinois at Chicago}
\author{David Bachman}
\author{Saul Schleimer}

\date{\today}

\begin{abstract}
  If a tangle, $K \subset \BB^3$, has no planar, meridional, essential
  surfaces in its exterior then thin position for $K$ has no thin
  levels.
\end{abstract}

\maketitle

\footnotetext[1]{To appear in the {\it Journal of Knot Theory and its Ramifications}.}

\section{Introduction}
In \cite{gabai:87}, D. Gabai introduced a complexity for embeddings of
knots in $S^3$, called {\em width}. A knot that has been isotoped to
have minimal possible width is said to be in {\em thin position}. In
\cite{thompson:97}, A. Thompson showed that if a knot, with no planar,
merdional essential surfaces in its exterior, is put in thin position
then all of its maxima appear above all of its minima. This
establishes a connection between {\em width}, and {\em bridge number},
a classical complexity introduced by H. Schubert \cite{schubert:54}. 

This paper gives a similar result for tangles in $\BB^3$. If a tangle
in $\BB^3$, with no planar, meridional, essential surfaces in its
exterior, is put in thin position then all of its maxima (if there are
any) appear above all of its minima. This is equivalent to saying that
thin position for such a tangle has no {\em thin levels} (see the next
section for the relevant definitions). This result has already been obtained when the
tangle is the one-skeleton of a triangulation (see \cite{schleimer:00}). 

The proof breaks into two cases. The first is handled by Thompson's
argument identically as in the situation of a knot. The second is
reminiscent of Claim~4.5 of Thompson's~\cite{thompson:94}. In this
case the {\em lightbulb trick} produces a contradiction to thin
position (as opposed to Thompson's ``{\em fluorescent lightbulb
trick}".) This change in technique is unavoidable; in
\cite{thompson:94} Thompson begins with a {\em thick level} of a
tangle, whereas we must start with a {\em thin level}. 

It has come to our attention that the second case may also be handled with an argument of S. Matveev's \cite{matveev:95}.

\section{Definitions}

Suppose that $K \subset \BB^3$ is a {\em tangle}, {\it i.e.} a properly embedded
1-manifold. Let $h: \BB^3 \rightarrow I$ be the standard height
function, so that $h^{-1}(t)$ is a 2-sphere for $t \in (0,1)$,
$h^{-1}(0)$ is a single point, and $h^{-1}(1) = \bdy \BB^3$.

Suppose that $h|K$, $h$ restricted to $K$, is a Morse function.  Let
$\{ c_i \}$ denote the critical values of $h|K$, and let $r_i$ be a
regular value immediately above $c_i$.

\begin{define}
The {\em width} of $K$ is the quantity 
\[w(K) = \sum_i |h^{-1}(r_i) \cap K|\] 
\end{define}

\begin{define}
The tangle $K$ is in {\em thin position} if $w(K) \leq w(K')$ for
any tangle $K'$ properly isotopic to $K$.
\end{define}

\begin{define}
A level, $L_i=h^{-1}(r_i)$, is {\em thin} if $|K \cap L_i| < |K \cap L_{i-1}|$ and $|K \cap L_i|<|K \cap L_{i+1}|$.
\end{define}

Let $M_K$ be the {\em exterior} of $K$, that is, $\BB^3$ with a small
neighborhood of $K$ removed. A surface properly embedded in $M_K$ is 
{\em essential} if it is incompressible and not boundary
parallel. Finally, a properly embedded surface is {\em meridional} if 
all boundary components also bound small meridional disks for $K$.

\section{The Main Theorem}

\begin{theorem}
\label{MainTheorem}
Let $K$ be a tangle in thin position in $\BB^3$. Suppose that $K$ has
a thin level. Then the exterior of $K$ contains a meridional, planar,
essential surface.
\end{theorem}

\begin{proof}
Pick a thin level, $S$. If $S$ does not intersect $K$ then $K$ is split, and $S$ is essential. If not, then note that $S_K = S \cap M_K$ is not boundary parallel in $M_K$. Compress $S$ as much as possible in the complement
of $K$ to obtain a collection of spheres, $P$, in $\BB^3$. For a
contradiction assume that each component of $P_K = P \cap M_K$ is
boundary parallel in $M_K$.

Suppose some component of $P_K$ is an annulus, parallel into the boundary
of a neighborhood of $K$. Pick an innermost such annulus, $P'_K$. Following Thompson's
proof of the main result of~\cite{thompson:97}, let $D$ be a disk in $\BB^3$ such that $\partial D=\delta \cup \gamma$, where $D \cap K=\delta$ and $D \cap P=D \cap P'=\gamma$. 
 
We now reverse the compressions used to obtain $P$ from $S$. Each time a compression is reversed, we attach a tube to some components of $P$. These tubes may intersect $D$, but only in its interior. In the end, $D$ persists as a disk in $\BB^3$ such that a collar of $\gamma$ in $D$ is disjoint from $S$. Isotope $\delta$ along $D$ to $\gamma$. This reduces width, a contradiction.

Suppose now there are no such annuli among the $n$ components of $P_K$, and each
component is a boundary parallel surface in $M_K$. Recall that
all boundary components of $P_K$ are meridional. Then each component
of $P$ is parallel into $\bdy \BB^3$ via an isotopy fixing $K$ setwise. This situation is
similar to Claim~4.5 from \cite{thompson:94}.

To fix notation, let $\{D_i\}_{i = 1}^n$ be a sequence of disks, and
$\{S_i\}_{i = 0}^n$ the sequence of surfaces, so that $S_0 = S$, $S_n = P$, and $S_i$ is obtained from $S_{i-1}$ by compressing along $D_i$ in the complement of $K$.  That is, remove a small neighborhood, $A_i$,
of $\bdy D_i$ from $S_{i-1}$. Construct $S_i$ by gluing two parallel
copies of $D_i$ onto $\bdy A_i$. Denote these by $B_i$ and $C_i$. So
$A_n \cup B_n \cup C_n$ bounds a ball homeomorphic to $D_n \times I$.
Finally, let $\alpha$ be the image of $\{\pt\} \times I \subset D_n  \times I$ by such a homeomorphism. The arc $\alpha$ is a {\em co-core} for the compression $D_n$.

Let $Q$ be the component of $S_{n-1}$ which meets $D_n$.  Let
$\{P_j\}_{j = 0}^{n}$ denote the components of $P = S_n$, numbered
consecutively, so that $P_0$ is innermost. Then there is a $j$ such
that compressing $Q$ along $D_n$ yields $P_j$ and $P_{j+1}$.  Let $N$
denote the submanifold of $\BB^3$ bounded by $P_j$ and $P_{j+1}$. Note
$\alpha \subset N$.

Choose a homeomorphism $g:S^2 \times I \rightarrow N$ such that $K \cap N$ is a collection of {\em straight arcs} and $\pi(g^{-1}(B_n)) = \pi(g^{-1}(C_n))$. Here a straight arc is one of the form $g(\{\pt\} \times I)$ and $\pi:S^2 \times I \rightarrow S^2$ is projection onto the first factor.

\begin{claim}
There is an isotopy, $\bar H: \BB^3 \times I \rightarrow \BB^3$,
fixing the complement of $N$ pointwise such that $\bar H(\alpha,1)$
is a straight arc.
\end{claim}

This is the usual {\em lightbulb trick}. See \cite{rolfsen:90}, for
example.  Since width is measured with respect to the height function,
$h$, we must use the reverse of $\bar H$ to move $K$. Let $H(\cdot,t) = \left(\bar H(\cdot,t)\right)^{-1}$. The isotopy $H$ moves $K$ to $H(K,1)$; this probably increases the width enormously. Another isotopy of $K$ is
needed to bring the width back under control.

Let $R = Q \cap S$ ($ = Q \setminus (\bigcup_{i=0}^{n-1} B_i \cup C_i)$) and set $R' = \bar H(R,1)$. Let $\{\beta _i\}_{i=1}^{m}$ denote the arcs of $K \cap N$.

\begin{claim}
There is an isotopy, $J: \BB^3 \times I \rightarrow \BB^3$, such that 
\begin{itemize}
    \item $J$ fixes the complement of $N$ pointwise,
      
    \item the arcs $J(\beta _i,1)$ lie within a very small
      neighborhood of $R'$, and
      
    \item the arcs $H(J(\beta _i,1),1)$ have only one critical point
      with respect to the height function $h$.
\end{itemize}

\end{claim}

        \begin{figure}[htbp]
        \psfrag{R}{$R'$}
        \psfrag{K}{$K$}
        \psfrag{J}{$J$}
        \vspace{0 in}
        \begin{center}
        \epsfxsize=3 in
        \epsfbox{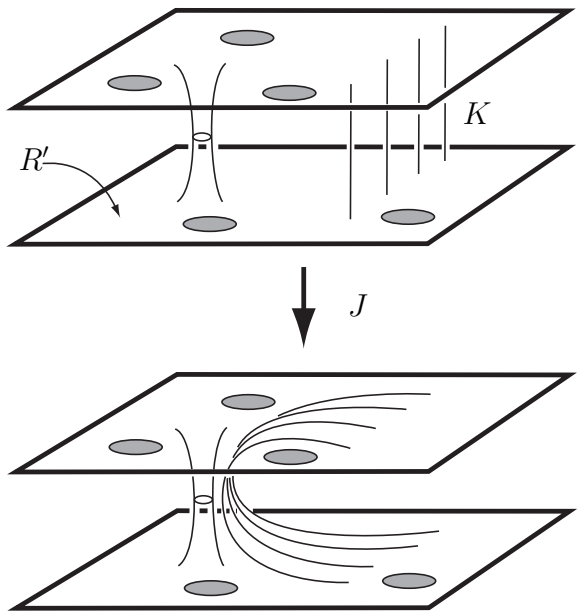}
        \caption{The isotopy, $J$, inside $N$.}
        \label{f:LightbulbIsotopy}
        \end{center}
        \end{figure}

\begin{proof}
Note that $\beta _i$ and $\bar H(\alpha,1)$ are straight. Hence,
there is an isotopy, $J$, supported in $N$, which moves the arcs
$\beta _i$ near $R'$ (see Figure \ref{f:LightbulbIsotopy}). The
image of $R'$ under $H(\cdot,1)$ is the surface, $R$, which is part
of the level, $S$. The foliation defined by the levels of $h$ is a
product near $R$, so the arcs $H(J(\beta _i,1),1)$ lie in a product
region. Adjust $J$ so that each arc has only one critical point with
respect to $h$.
\end{proof}

For each $i \in \{1, \ldots, m\}$, let $\beta'_i = H(J(\beta_i,1),1)$. For each
$j \in \{1, \ldots, m\}$, let $K_j$ be the 1-manifold which consists
of the subset of $K$ which lies outside of $N$, together with the arcs
$\{\beta_i'\}_{i=1}^{j}$ and the arcs $\{\beta _i\}_{i=j+1}^{m}$.
Finally, let $K_0=K$. Hence, for each $i \in \{1, \ldots, m\}$, $K_i$
can be obtained from $K_{i-1}$ by isotoping the arc, $\beta _i$, to be
near $S$ (refer to Figure \ref{f:knot2}.) $K_i$ is not necessarily
isotopic to $K_{i-1}$. However, $K_m$ is isotopic to $K_0$, since
$K_m=H(J(K,1),1)$.

We now assume that the arcs $\beta_i$ lie above $S$ with respect to the height function $h$. (The other case is identical.)

\begin{claim}
  For each $i \in \{1, \ldots, m\}$, $w(K_{i-1}) \geq w(K_i)$. Here
  equality implies that the first critical point of $h|K_{i-1}$ above
  $S$ is a maximum.
\end{claim}

Proving this claim will produce the desired contradiction. For if
$w(K_0) > w(K_m)$ then $K=K_0$ was not thin. On the other hand, if
$w(K_0) = w(K_m)$ then $S$ was not a thin level.

\begin{proof}
Let $\{r_j\}$ denote regular values of $h|K_i$ which occur just above
each critical value. Let $b'$ be the unique critical value of $\beta' _i$ and let $b$ be the largest critical value of $\beta_i$.  Suppose $r_l$ is just above $b'$ and $r_{l+k}$ is the highest regular value
below $b$. Finally, let $s$ be a regular value just above $b$ (see schematic Figure \ref{f:knot2}).

        \begin{figure}[htbp]
        \psfrag{c}{$b'$}
        \psfrag{b}{$b$}
        \psfrag{s}{$s$}
        \psfrag{S}{$S$}
        \psfrag{1}{$r_l$}
        \psfrag{2}{$r_{l+1}$}
        \psfrag{k}{$r_{l+k}$}
        \psfrag{K}{$K_{i-1}$}
        \psfrag{L}{$K_i$}
        \psfrag{B}{$\partial B$}
        \psfrag{e}{$\beta _i$}
        \psfrag{f}{$\beta '_i$}
        \vspace{0 in}
        \begin{center}
        \epsfxsize=4 in
        \epsfbox{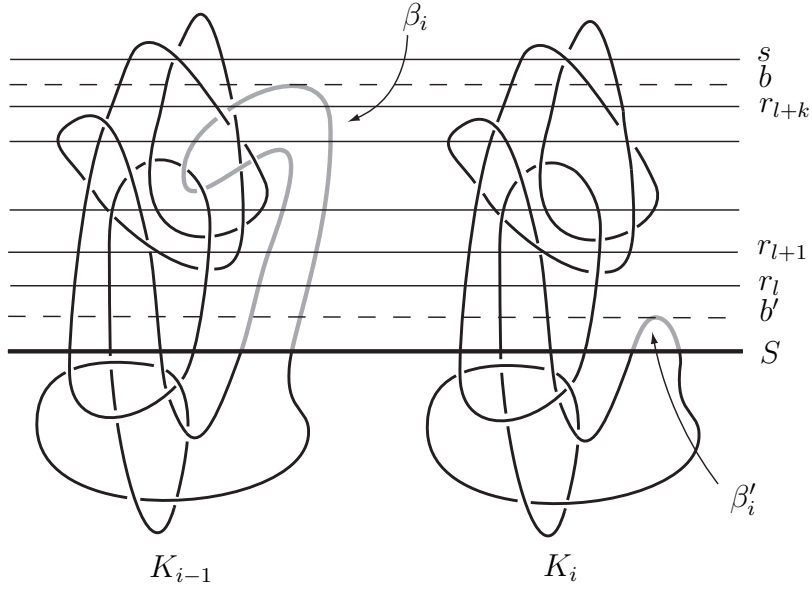}
        \caption{Comparing widths.}
        \label{f:knot2}
        \end{center}
\end{figure}

The width of $K_i$ is then $\sum_j |h^{-1}(r_j) \cap K_i|$ and the width of $K_{i-1}$ is at least
\[\left( \sum \limits _j |h^{-1}(r_j) \cap K_{i-1}|\right) + |h^{-1}(s) \cap K_{i-1}| - |h^{-1}(r_l) \cap K_{i-1}|.\]
Note that $|h^{-1}(r_j) \cap K_{i-1}| = |h^{-1}(r_j) \cap K_i|$ for all values of $j > l+k$ or $j < l$. Hence the difference  $w(K_{i-1}) - w(K_i)$ is at least:

\hspace{1cm}\lefteqn{\sum \limits _{j=l+1}^{l+k} \left( |h^{-1}(r_j) \cap K_{i-1}| - |h^{-1}(r_j) \cap K_i| \right)}
\vspace{-0.25cm}
\[+ |h^{-1}(s) \cap K_{i-1}| - |h^{-1}(r_l) \cap K_{i-1}|\]
\vspace{-0.33cm}

Also, for each $j \in \{l$+$1,...,l$+$k\}$ we have $|h^{-1}(r_j) \cap K_{i-1}| \geq 2 + |h^{-1}(r_j) \cap K_i|$, so the above difference in width is at least

\vspace{-0.33cm}
\[\Delta=2k + |h^{-1}(s) \cap K_{i-1}| - |h^{-1}(r_l) \cap K_{i-1}|.\]

Note that $\Delta$ is equal to zero only when all of the critical values of $h|K_{i-1}$ in the interval $[b',b]$ are maxima. Otherwise, $\Delta$ is strictly greater than zero. In particular, if $w(K_{i-1}) - w(K_i) = 0$ then the first critical point of $h|K_{i-1}$ above $S$ is a maximum.
\end{proof}

This completes the proof of Theorem~\ref{MainTheorem}.
\end{proof}

\bibliographystyle{plain}

\bibliography{thinball}

\end{document}